\documentclass[11pt]{amsart}
\usepackage[margin=2.5cm]{geometry}

\usepackage[english]{babel}
\usepackage{amsthm,amsmath,amsfonts,amssymb,commath}   
\usepackage{mathrsfs}
\usepackage{upgreek}
\usepackage{thmtools}
\usepackage{lipsum}

\newdimen\parindentt
\parindentt=\parindent
\advance\parindentt by 5mm

\usepackage{geometry}
\usepackage{indentfirst}   

\usepackage[usenames,svgnames,dvipsnames]{xcolor}
\definecolor{verde}{rgb}{0,0.5,0}

\definecolor{laranja}{rgb}{0.95,0.45,0}

\definecolor{vermelho}{rgb}{0.666,0,0}

\usepackage[all]{xy}
\usepackage[utf8]{inputenc}
\usepackage{mathptmx}
\usepackage{mathrsfs}
\usepackage{amsaddr}

\usepackage{hyperref,cleveref}
\usepackage{enumerate,multicol}
\usepackage{tikz-cd}
\usepackage{etoolbox}

\newtheorem{theorem}{Theorem}[section]

\newtheorem{corollary}[theorem]{Corollary}

\theoremstyle{definition}
\newtheorem{definition}[theorem]{Definition}
\newtheorem{example}[theorem]{Example}

\theoremstyle{remark}
  \newtheorem{remark}[theorem]{Remark}
  \AtBeginEnvironment{remark}{%
  \pushQED{\qed}%
}
\AtEndEnvironment{remark}{\popQED\endremark}

\numberwithin{equation}{section}

\DeclareMathOperator{\ga}{\mathsf{g}}
\DeclareMathOperator{\aga}{\mathsf{h}}

\def\cal#1{\mathcal{#1}}
\def\bb#1{\mathbb{#1}}
\def\lie#1{\mathfrak{#1}}

\usepackage[all]{xy}

\newcommand{\h}{\frac{1}{2}}

\newcommand\xqed[1]{%
  \leavevmode\unskip\penalty9999 \hbox{}\nobreak\hfill
  \quad\hbox{#1}}
\newcommand\demo{\xqed{$\blacksquare$}}

\tikzset{
  symbol/.style={
    draw=none,
    every to/.append style={
      edge node={node [sloped, allow upside down, auto=false]{$#1$}}}
  }
}

\usepackage[mathlines]{lineno}

\title{A note on the vertizontal curvature of fat bundles}
\author{Leonardo F. Cavenaghi}
\address{Instituto de Matemática, Estatística e Computação Científica -- Unicamp, Rua Sérgio Buarque de Holanda, 651, 13083-859, Campinas, SP, Brazil}
\email{leonardofcavenaghi@gmail.com}

\begin{document}
	

\begin{abstract}
In \cite{Ziller_fatnessrevisited}, W. Ziller poses a compelling question: given a fat principal $G$-bundle $(P, \ga) \rightarrow (B, \aga)$ with $\dim G = 3$, and $\ga$ representing a Riemannian submersion metric ensuring that the $G$-orbits are totally geodesic, can one modify $\aga$ to render all vertical curvatures equal to $1$? In this note, we establish a rigidity result for fat Riemannian foliations with bounded holonomy and a specific curvature constraint. Our result addresses Ziller's question for fat fiber bundles with compact structure groups, considering connected compact total spaces under a curvature constraint that holds on various examples, such as locally symmetric spaces. Additionally, we assume that all vertizontal curvatures coincide at a point. 
\end{abstract}
	
	\maketitle
	
	\section{Introduction and Preliminaries}

Let $\pi: F\hookrightarrow (M,\ga) \rightarrow (B,\aga)$ be a Riemannian submersion where $F$ represents the fiber and $B$ is the base. Let $\cal V$ denote the vertical bundle of $\pi$, which contains vectors tangent to the fibers. We term the \emph{horizontal bundle} its $\ga$-orthogonal complementary bundle, denoted as $\cal H$ and isometric to $(TB,\aga)$. We say that the Riemannian submersion $\pi$ is \emph{fat} if, for every point $x$ in $M$ and every nonzero vector $X$ in $\cal H_x$, the following condition holds for a local horizontal extension $\widetilde X$ of $X$:
\begin{equation}\label{eq:fat}
[\widetilde X,\cal H_x]^{\mathcal V} = \cal V_x.
\end{equation}
The left-hand side in Equation \eqref{eq:fat} is $2A_X\cal H_x$ where $A:\cal H_x\times \cal H_x\rightarrow \cal V_x$ stands for the O'Neill tensor of the submersion $\pi$.

Initially introduced in \cite{weinstein1980fat}, fat bundles offer a valuable framework for exploring the conditions necessary for metrics with positive sectional curvature on the total space of fiber bundles, particularly those with totally geodesic fibers. According to Gray \cite{gray1967pseudo} or O'Neill \cite{oneill}, in the presence of totally geodesic fibers, the unreduced sectional curvature of $\ga$ at any nontrivial plane $X\wedge V$ for $X\in \mathcal{H}, V\in \mathcal{V}$ is given by
\[K_{\ga}(X,V)=|A^*_XV|_{\ga}^2,\]
where $A^*_X$ is the $\ga$-dual to $A_X$. This underscores the necessity of the fatness condition for a metric with positive sectional curvature and totally geodesic fibers on $M$. A plane of the form $\sigma = X\wedge V$ is termed a \emph{vertizontal plane}. In this manuscript, we adopt the following convention for the Riemannian tensor
\[R_{\ga}(X,Y)Z=\nabla_X\nabla_YZ-\nabla_Y\nabla_XZ-\nabla_{[X,Y]}Z,~\forall X, Y, Z\in TM.\]
Therefore, we define the sectional curvature of a plane $X\wedge Y$ to be 
\[\mathrm{sec}_{\ga}(X,Y)=R_{\ga}(X,Y,Y,X)/|X|^2|Y|^2-\langle X,Y\rangle^2\]
where $R_{\ga}(X,Y,Y,X)=\ga\left(R_{\ga}(X,Y)Y,X\right)$. We sometimes employ the notation $\langle\cdot,\cdot\rangle$ to shortly denote $\ga(\cdot,\cdot)$.

As highlighted in \cite[Proposition 2.5, p.8]{Ziller_fatnessrevisited}, for principal fat bundles with fiber $F = \mathrm{S}^3$ or $\mathrm{SO}(3)$, the base manifold $B$ must have a dimension multiple of $4$. Furthermore, the only known instances of principal fat bundles with these fibers occur in $3$-Sasakian manifolds. A manifold $M$ is labeled $3$-Sasakian if it admits an isometric almost free action by $\mathrm{S}^3$ or $\mathrm{SO}(3)$ with totally geodesic orbits of curvature $1$, satisfying $R_{\ga}(U^*, X)Y = -\ga(U^*, Y)X{+}\ga(X, Y)U^*$ for all $U \in \mathfrak{g}$, the Lie algebra of $G$. The superscript $*$ denotes the action vector field derived from $U$. For these metrics, all vertizontal curvatures are equal to $1$. In W. Ziller's notes on fat bundles, it is asked: ``Given a fat principal $G$-bundle with $\dim G = 3$, is it possible to change the metric on $B$ such that all vertizontal curvatures are equal to $1$?'' This note shows that this assertion holds for fat fiber bundles with compact structure groups that are locally symmetric spaces if all vertizontal curvatures equal $1$ at a point. This follows from the more general rigidity result, which states for compact connected Riemannian manifolds with fat Riemannian foliations with bounded holonomy (Theorem \ref{thm:main}) if at a point the verzitonal curvatures coincide, it coincide on the whole manifold under an additional curvature constraint. To make a better sense of this, we introduce further notation.

The connected components of the fibers of the Riemannian submersion $\pi:(M,\ga)\rightarrow (B,\aga)$ consist of an example of a Riemannian foliation. However, a Riemannian foliation does not always need to be given in this form. In the following, we employ the same terminology used in Riemannian submersions to Riemannian foliations. Namely, tangent vectors to leaves are named vertical. For any $x\in M$, the vector space collecting vertical vectors is named the vertical space and denoted by $\cal V_x$. Its $\ga$-complementary vector space is denoted as $\cal H_x$ and named the horizontal space. Its elements are named horizontal vectors. We also maintain the notation $A_XY=\h [X,Y]^{\cal V}$ for the O'Neill tensor, where $X,Y\in \cal H_x$ and $A^*_X$ its $\ga$-dual once $X\in \cal H_x$ is fixed. Since fatness is solely related to a choice of complementary distribution to the fibers (in the Riemannian submersion case), we employ the same definition for Riemannian foliations. To wit, 

\begin{definition}
    Let $(M,\ga,\cal F)$ be a Riemannian manifold with a Riemannian foliation $\cal F$. We say that $\cal F$ is fat if, for every $x\in M$ and every nonzero vector $X$ in $\cal H_x$, the following condition holds for a local horizontal extension $\widetilde X$ of $X$:
\begin{equation}
[\widetilde X,\cal H_x]^{\mathcal V} = \cal V_x.
\end{equation}
\end{definition}

\begin{definition}[Holonomy Fields and Dual Holonomy Fields]
Let $(M,\ga,\cal F)$ be a compact Riemannian manifold with a Riemannian foliation $\cal F$. Let $c$ be a horizontal curve on $M$. The \emph{holonomy field} $\xi$ along $c$ is a vertical field satisfying
\begin{equation}\label{eq:holonomyfield}
\nabla_{\dot{c}}\xi=-A^*_{\dot{c}}\xi-S_{\dot c}\xi,
\end{equation}
where $S$ is the shape operator, i.e., the $\ga$-dual to the second fundamental form of the leaves. Likewise, a vertical field $\nu$ along $c$ is called a \emph{dual holonomy field} if
\begin{equation}\label{eq:dualhol}
\nabla_{\dot{c}}\nu=-A^*_{\dot{c}}\nu+S_{\dot c}\nu.
\end{equation}
\end{definition}
 
 For a horizontal curve $c:[0,1]\to M$, let $h:\cal V_{c(0)}\to \cal V_{c(1)}$ be the linear isomorphism given by $h(\xi_0)=\xi(1)$, where $\xi(t)$ is the holonomy field along $c$ with initial condition $\xi(0)=\xi_0$. Following \cite{speranca2017on}, we name $h$ an \textit{infinitesimal holonomy transformation}.
\begin{definition}[Bounded Holonomy]\label{def:bhol}
We say that a Riemannian foliation has bounded holonomy if there is a constant $L$ such that, for every holonomy field $\xi$, $|\xi(1)|\leq L|\xi(0)|$.
\end{definition}
\begin{example}[Proposition 3.4 in \cite{speranca2017on}]
     If the structure group of the Riemannian submersion $\pi:(M,\ga)\rightarrow (B,\aga)$ with compact connected total space is compact, the Riemannian foliation induced by the connected component of the fibers has bounded holonomy.
\end{example}

Let $(M,\ga,\cal F)$ be a compact connected Riemannian manifold with a Riemannian foliation $\cal F$ with totally geodesic leaves. A triple $\{X,V,\cal A\}$ is said to be a \emph{good triple} if $\exp(tV (s)) = \exp(sX(t))$ for all $s, t \in \bb R$, where $V(s), X(t)$ denote the Jacobi fields along $\exp(sV)$ and $\exp(tX)$, respectively, that satisfy $V(0) = V , X(0) = X$ and $V'(0) = \cal A = X'(0)$. Proposition 1.4 in \cite{tappmunteanu2} ensures that since the leaves of $\cal F$ are totally geodesic, then the triple $\{X, V,-A^*_XV\}$ consists of a good triple for every $X\in \cal H_x$ and $V\in \cal V_x$. For our purposes, we need the following:
\begin{theorem}[Munteanu--Tapp, \cite{tappmunteanu2}]
On any locally symmetric space $(M,\ga)$ it holds that $R_{\ga}(X,V)(A^*_XV) = 0$ where $R_{\ga}$ is the Riemannian curvature tensor of $\ga$.
\end{theorem}

We are in a position to state our results:
\begin{theorem}\label{thm:main}
Let $(M,\ga,\cal F)$ be a compact connected Riemannian manifold with a fat Riemannian foliation $\cal F$ with totally geodesic leaves and bounded holonomy.  If for every good triple $\{X,V,-A^*_XV\}$ we have $R_{\ga}(X, V)A^*_XV=0$ then the vertizontal curvatures of $\ga$ are determined at a point. This holds, for instance, if $(M,\ga)$ is isometric to a locally symmetric space.

If, in addition, the vertizontal curvatures of $\ga$ coincide at a point, this holds for every point. In particular, a scaling of $\ga$ exists for which all vertizontal curvatures at every point equal $1$. 
\end{theorem}
As an immediate corollary we obtain a related result to W. Ziller's question.
\begin{corollary}
    Let $(M,\ga)$ be a compact connected Riemannian manifold with a free isometric action by $G=\mathrm{S}^3,~\mathrm{SO}(3)$. Assume that the orbits of $G$ in $(M,\ga)$ are totally geodesic and pattern a fat Riemannian foliation in $(M,\ga)$. If
    \begin{enumerate}[(a)]
        \item for all good triples $\{X,V,-A^*_XV\}$ we have $R_{\ga}(X,V)A^*_XV=0$
        \item there exists a point $x\in M$ where all the vertizontal curvatures are equal to $1$
    \end{enumerate}
    then $(M,\ga)$ is a $3$-Sasakian manifold.
\end{corollary}

Applying the former to the examples discussed in \cite{Ziller_fatnessrevisited}, one recovers the following by observing that the total spaces below are locally symmetric spaces with a point with positive vertizontal curvature equal to $1$.
\begin{corollary}\label{cor:ziller}
    In the following fat $\mathrm{S}^3,~\mathrm{SO}(3)$-principal bundles a Riemannian submersion metric exists with constant vertizontal curvature equal to $1$. In particular, the total spaces are $3$-Sasakian manifolds.
    \begin{enumerate}[(a)]
    \item $\mathrm{S}^3\rightarrow \mathrm{S}^{4n+3}\rightarrow \bb HP^n$
        \item $\mathrm{SO}(3)\rightarrow T_1\mathbb{C}P^n\rightarrow \mathrm{Gr}_2(\bb C^{n+1}),~n\geq 2$
        \item $\mathrm{SO}(3)\rightarrow \mathrm{SO}(n)/\mathrm{S}^3\times\mathrm{SO}(n-4)\rightarrow \Lambda_{4,n}^{\circ},~n\geq 5$ (the oriented $4$-planes in $\bb R^n$)
        \item $\mathrm{S}^3\rightarrow \mathrm{Sp}(n+1)/\mathrm{Sp}(n)\rightarrow \bb HP^n$
        \item $\mathrm{SO}(3)\rightarrow \mathrm{G}_2/\mathrm{S}^3_{\pm}\rightarrow \mathrm{G}_2/\mathrm{SO}(4)$ where $\mathrm{S}^3_{\pm}\subset \mathrm{SO}(4)$ are the two simple factors.
    \end{enumerate}
\end{corollary}

	\section{Proof of Theorem \ref{thm:main}}

	Let $(M,\ga)$ be a compact connected Riemannian manifold with a fat Riemannian foliation $\cal F$ with totally geodesic leaves and bounded holonomy. Theorem 6.2 in \cite{speranca2017on} implies that any two points in $(M,\ga)$ can be joined by a  smooth horizontal geodesic.
 
 Take any two points $x\neq y\in M$ and let $c : [0,d_{\ga}(x,y)]\rightarrow M$ to be a smooth horizontal geodesic joining $x$ to $y$. Varying every possible $y$ once fixed $x$ (or \emph{vice-versa}) collects every possible horizontal direction. Let $V$ be a vertical vector at $x$. \cite[Proposition 1.4]{tappmunteanu2} implies that for $X=\dot{c}(0)$ the set $\{X,V,-A^*_XV\}$ is a good triple, the same holding for the set $\{\dot{c}(t),\nu(t),-A^*_{\dot{c}(t)}\nu(t)\}$ for every $t\in [0,d_{\ga}(x,y)]$ where $\nu(t)$ as a dual-holonomy field extending $V$ along $c$. 
 
 Following \cite[Theorem 2]{oneill} we have
 \begin{align*}
    R_{\ga}(\dot{c},A^*_{\dot{c}}\nu,\dot{c},\nu) &= -\langle(\nabla_{\dot{c}}^{\cal V}A)_{\dot{c}}A^*_{\dot{c}}\nu,\nu\rangle\\
    &=-\langle \nabla_{\dot{c}}(A_{\dot{c}}A^*_{\dot{c}}\nu)-A_{\dot{c}}(\nabla_{\dot{c}}A^*_{\dot c}\nu),\nu\rangle\\
    &=-\langle \nabla_{\dot{c}}(A_{\dot{c}}A^*_{\dot{c}}\nu),\nu\rangle+\langle \nabla_{\dot{c}}A^*_{\dot{c}}\nu,A^*_{\dot{c}}\nu\rangle\\
    &=-\frac{d}{dt}\langle A_{\dot{c}}A^*_{\dot{c}}\nu,\nu\rangle+\langle A^*_{\dot{c}}\nu,A^*_{\dot{c}}(\nabla_{\dot{c}}\nu)\rangle+\langle \nabla_{\dot{c}}A^*_{\dot{c}}\nu,A^*_{\dot{c}}\nu\rangle\\
    &=-\frac{d}{dt}|A^*_{\dot{c}}\nu|^2+\langle \nabla_{\dot{c}}A^*_{\dot{c}}\nu,A^*_{\dot{c}}\nu\rangle.
 \end{align*}
According to our hypotheses we have $R_{\ga}(\dot{c}, A^*_{\dot{c}}\nu,\nu,\dot{c})=0$. Thus, the former sequence of equalities ensures that
\[\frac{d}{dt}|A^*_{\dot{c}}\nu|^2=\langle \nabla_{\dot{c}}A^*_{\dot{c}}\nu,A^*_{\dot{c}}\nu\rangle.\]
However,
\begin{align*}
	 \frac{d}{dt}K_{\ga}(\dot{c},\nu) = \frac{d}{dt}|A^*_{\dot{c}}\nu|_{\ga}^2 &= 2\langle A^*_{\dot{c}}\nu,\nabla_{\dot{c}}(A^*_{\dot{c}}\nu)\rangle.
	\end{align*}
 Therefore, $\frac{d}{dt}K_{\ga}(\dot{c},\nu) =0$, this proves the first assertion.
 
 To the last assertion of the theorem, observe that if there is a point $x\in M$ where all vertizontal curvature coincides, the former computation ensures that for any other $y$, these curvatures coincide with that at $x$. More precisely, assuming that at $x$ we have $K_{\ga}(X,V)=1$ for every non-degenerate plane $X\wedge V$, the set $\cal M:=\{x\in M:K_{\ga}(X,V)=1~\forall X\in \cal H_x,~V\in \cal V_x,~X\wedge V\neq 0,~\text{such that}~K_{\ga}(X,V)=1\}=M$. Indeed, we show that $\cal M$ is open and closed. Since it is non-empty the connectedness of $M$ will ensure the result.

 Since horizontal and vertical distributions globally form sub-bundles of $TM$, it is obvious from the continuity of $K_{\ga}$ in the two-plane Grasmannian that $\cal M$ is closed. It is left to prove it is open. However, this already follows from the former step. Due to Theorem 6.2 in \cite{speranca2017on} one has that a geodesic ball centered at $x$ contains a ball $\cal B$ in $M$ where two points can be joined by a horizontal smooth geodesic. Thus, according the first part of the proof, for any $X\in \cal H_x$ and any $V\in \cal V_x$, the curvature $K_{\ga}(\dot c,\nu)$ is constant along the horizontal geodesic generated by $X$ where $\nu$ is a dual-holonomy field extension of $V$ along $c$. Therefore, $\cal M$ contains $\cal B$. Thus, $\cal M$ is open, as desired.
 
 \demo

 \section{Proof of Corollary \ref{cor:ziller}}
 We recall that a homogeneous bundle is of the form
 \[K/H\rightarrow G/H\rightarrow G/K\]
where $H<K<G$. If $H$ is normal in $K$, then the former homogeneous bundle is said to be \emph{principal homogeneous}. 

Once chosen an $\mathrm{Ad}(K)$-invariant splitting $\lie g = \lie k\oplus \lie m$, there exists an invariant connection $\theta$ on $G\rightarrow G/K$ that induces a horizontal homogeneous distribution $\widetilde{\cal H}$ in $G/H\rightarrow G/K$ satisfying $\widetilde{\cal H}_{(H)}\cong \lie m$. Moreover, we can induce homogeneous metrics on $K/H$ and $G/K$ via the identifications $T_{(H)}K/H\cong \lie p,~T_{(K)}G/K\cong \lie m$ after choosing an $\mathrm{Ad}(H)$-invariant splitting $\lie k=\lie h\oplus \lie p$ and $\mathrm{Ad}(H),~\mathrm{Ad}(K)$-invariant metrics on $\lie p$ and $\lie m$, respectively. Gathering this all, we can define a mertic $\ga$ on $G/H$ satisfying $T_{(H)}G/H\cong \lie p\oplus \lie m$. As Proposition 3.1 in \cite{Ziller_fatnessrevisited} ensures, we have $K_{\ga}(X,V)=|[X,V]|_{\ga}^2$ for $X\in \lie m,~V\in \lie p$.

In the case $K/H=\mathrm{S}^3,~\mathrm{SO}(3)$, the former homogeneous fibrations can be described in terms of homogeneous quaternionic symmetric spaces $G/K.\mathrm{Sp}(1)$ (named Wolf spaces). Proposition 3.4 in \cite{Ziller_fatnessrevisited} states that if $G/K.\mathrm{Sp}(1)$ is a quaternionic symmetric space, then the principal bundle fibration $K.\mathrm{Sp}(1)/K\rightarrow G/K\rightarrow G/K.\mathrm{Sp}(1)$ is fat. Furthermore, that $K.\mathrm{Sp}(1)/K=\mathrm{S}^3$ or $\mathrm{SO}(3)$ and $\mathrm{Sp}(1)/K$ is 3-Sasakian. Every listed example in the statement of Corollary \ref{cor:ziller} is of this form. It is also straightforward to check that the corresponding total spaces are realized as symmetric spaces; see pp.18 in \cite{Ziller_fatnessrevisited}.

 Let $\bb H$ stands to the quaternions and consider the Hopf fibration $\mathrm{Sp}(1)\rightarrow \mathrm{S}^{4n+1}\rightarrow \bb{H}^n$. Adapting the metric construction above to this fibration one trivially has that the vertizontal curvatres for this fibration are equal to $1$. This fact concludes the proof since \cite[pp. 17]{Ziller_fatnessrevisited} ensures that the restriction of the isotropy representation of $K.\mathrm{Sp}(1)$ on $G/K.\mathrm{Sp}(1)$ to the subgroup $\mathrm{Sp}(1)$ is equivalent to the usual Hopf action of $\mathrm{Sp}(1)$ on $\bb{H}^{n}$. Thus, for the metric $\ga$ constructed as above we have that the vertizontal curvatures coincide with that of the quaternionic Hopf fibrations. \demo

   \section*{Ackowledgements}
   L.F.C. is funded by The São Paulo Research Foundation (FAPESP), grants 2022/09603-9 and 2023/14316-1.

\end{document}